\documentclass[preprint,12pt]{elsarticle}
\usepackage{subcaption}
\usepackage{amssymb}
\usepackage{array}
\usepackage{graphicx}
\usepackage{float}
\usepackage{subfig}
\usepackage{longtable}
\usepackage{algorithmic}
\usepackage{algorithm}
\usepackage{amsmath}
\usepackage{amsthm}
\usepackage{caption}
\usepackage{url}
\usepackage{color}
\usepackage{booktabs}
\usepackage{bm}
\usepackage{tablefootnote}
\newcommand{\svhline}{\noalign{\hrule height 1.5pt}}
\theoremstyle{definition}
\newtheorem{remark}{Remark}
\journal{Elsevier}
\date{}
\begin{document}
\begin{frontmatter}
\title{Weighted Quadrature on Unstructured Splines} 
\author{Ji Sheng$^{a}$, Xiaodong Wei$^{a,\ast}$, Falai Chen$^{b}$} 
\affiliation[ a ]{organization={Shanghai Jiao Tong University Global College},
            addressline={Shanghai Jiao Tong University}, 
            city={Shanghai},
            postcode={200240}, 
            country={China}}
\affiliation[ b ]{organization={School of Mathematical Sciences},
            addressline={University of Science and Technology of China}, 
            city={Hefei},
            postcode={230026}, 
            country={China}}
\begin{abstract}
This work presents a weighted quadrature (WQ) method to fast assemble Galerkin matrices based on unstructured spline surfaces. The method is developed upon a particular variant of unstructured splines, namely the bicubic analysis-suitable unstructured T-splines (ASUTS). While existing WQ approaches have significant speedup for structured splines (e.g., B-splines), their extension to unstructured splines faces several  challenges: (1) lack of a global parametric domain for defining quadrature points, (2) a varying number of basis functions across elements that complicates the determination of the optimal number of quadrature points, and (3)  ill-conditioned underdetermined linear systems that must be solved to find the quadrature weights. To solve these issues, we first define the WQ rule directly
 in the physical domain. Second, we specify the number of quadrature points function-wise (rather than element-wise), which naturally satisfies the well-posedness condition, namely the number of unknown weights no less than that of exactness constraints. Third, we employ the truncated Singular Value Decomposition to improve the conditioning of the underdetermined systems by discarding extremely small singular values, which are caused by the splines around extraordinary points. Several different model problems are studied, such as Poisson's problem, the biharmonic problem, and the nonlinear heat transfer problem. In the end, a variety of numerical tests are performed to demonstrate the accuracy and efficiency of the proposed method.
\end{abstract}
\end{frontmatter}
\section{Introduction}
Isogeometric Analysis (IGA) was proposed to seamlessly integrate Compu\-ter-Aided Design (CAD) and Computer-Aided Engineering (CAE) \cite{HUGHES20054135}. It directly adopts spline-based geometric representations, such as non-uniform rational B-splines (NURBS) \cite{HUGHES20054135} and T-splines \cite{BAZILEVS2010229,SCOTT}, for both geometric modeling and numerical analysis. However, the high-order nature of splines poses a significant challenge to the efficient formation of Galerkin matrices \cite{Mant2014}. The conventional elementwise assembly
 combined with the standard Gaussian quadrature (GQ) leads to a
 high complexity of $\mathcal{O}(p^{3d})$ per element, where $d$ is the spatial dimension and $p$ is the spline degree.

Numerous studies have been dedicated to mitigating this computational bottleneck. For instance, generalized Gaussian quadrature was developed to extend classical Gaussian quadrature from polynomial spaces to more general function spaces such as splines \cite{GGQ1996,EQ2010,Optimal2016,Gaussian2016}. Reduced quadrature  aims to reduce the number of the standard Gaussian quadrature points while maintaining a sufficient level of integration accuracy \cite{SCHILLINGER20141,HIEMSTRA2017966}. Sum factorization fully exploits the tensor-product structure of B-splines/NURBS by factorizing common factors to eliminate redundant calculations \cite{2015Efficient}. In  \cite{lookup2015,MANTZAFLARIS20171062,PAN2020113005}, an integral  is approximated as a linear combination of integrals of B-spline (or related derivatives) products, which are precomputed and stored in a look-up table, with coefficients determined through interpolation. Weighted quadrature (WQ) constructs  quadrature  rules for each test function \cite{F2017Fast,Ren2019Fast,cal2019Fast}, so matrix assembly proceeds row-wise rather than element-wise, where sum factorization can also be leveraged to further enhance efficiency. Since it features low complexity of $\mathcal{O}(p^{d+1})$ per element and only solves local problems to determine the weights, WQ has been studied in a variety of contexts, such as boundary element methods \cite{AIMI2018327,Falini2019,Falini2019WQ}, hierarchical B-splines \cite{GIANNELLI2022115465,PAN2021113278}, and the immersed method \cite{MARUSSIG2023116397}. We also follow this direction for the same reason.

However, the vast majority of existing work on WQ has focused on B-splines/NURBS. They have a tensor-product structure and global parametric domains, so they can only model a limited family of geometries.  On the other hand, unstructured splines \cite{COLLIN201693,Wei2021,Wei2022AnalysissuitableUT,Casquero2020SeamlessIO} offer watertight geometric representations for complex geometries and thus serve as highly promising candidates for IGA to fundamentally integrate CAD with CAE. However, fast formation has not yet been studied in this context to the best of the authors' knowledge. The key characteristic of unstructured splines is the presence of \textit{extraordinary points} (EPs), which pose several fundamental barriers for the adoption of WQ.
 First, there are no global parametric domains in unstructured splines. Instead, parametric domains are local to elements. Second, the number of basis functions varies significantly across different types of elements, complicating the determination of the optimal number of quadrature points while satisfying the well-posedness condition, namely the number of unknown weights no less than that of the exactness constraints. Third, the underdetermined linear systems used to solve for the weights are ill-conditioned because of the unstructured splines around EPs. To address these issues, we propose a WQ method based on unstructured splines. In particular, we build our method upon \textit{bicubic} analysis-suitable unstructured T-splines (ASUTS) \cite{Wei2022AnalysissuitableUT,Casquero2020SeamlessIO}, a representative type of unstructured splines that possesses optimal convergence rates and mixed types of splines. Our choice of the bicubic case is motivated by the fact that it dominates geometric modeling of complex geometries in practice. Also note that while T-junctions are a primary feature of ASUTS, we focus only on EPs in this work.  Extensions to alternative spline constructions and polynomial degrees are possible.

The remainder of the paper is organized as follows. Section 2 reviews WQ on B-splines. Section 3 introduces the basics of ASUTS. Section 4 details the development of WQ on ASUTS. Section 5 introduces the extension to the nonlinear heat transfer problem. In Section 6, we demonstrate the effectiveness and efficiency of the proposed approach through a variety of numerical tests. Finally, Section 7 concludes the work and outlines promising future directions.

\section{Weighted quadrature on B-splines}
In this section, we introduce the core idea of weighted quadrature (WQ). For simplicity, we focus on mass matrices assembled with uniform univariate B-splines. WQ in the bivariate or trivariate case follows a tensor-product manner. More general matrices will be discussed in the subsequent sections.

B-spline basis functions $B_k(u)$ of degree $p$ are defined on a non-decreasing sequence of real numbers called a knot vector. For simplicity, we consider a uniform knot vector,
\begin{equation}
\{ 0, 1, 2, \ldots, m \} ,
\end{equation}
where $m>p+1$ is a given positive integer. Each B-spline $B_k$ has a local support over interval $(k, k+p+1)$, denoted as $\mathrm{supp}(B_k)$. The integral of interest is
\begin{equation}
m_{ij}=\int_0^mc(u)B_i(u)B_j(u)du,\quad i, j = 1,2,\ldots, n,
\end{equation}
where $c(u)$ contains the information of geometric mapping  (e.g.,  Jacobian determinant) and/or the involved partial differential equation (PDE), $n$ is the number of  B-splines, and $B_i$ and $B_j$ are the test and trial functions, respectively.

WQ  is defined for each test function $B_i$, which is eventually incorporated as part of quadrature weights.  $c(u)$ in Eq. (2) is assumed to be sufficiently smooth and thus neglected during the development of WQ. Specifically, for test function $B_i$, we have
\begin{equation}
\int_0^mB_j(u)\bigg(B_i(u)du\bigg)\approx\mathbb{Q}_i\bigg(B_j(u)\bigg)=\sum_qw_{i,q}B_j(u_{i,q}),\,\,\,\forall j,
\end{equation}
where $\{u_{i,q}\}$ is a predetermined set of quadrature points and $\{w_{i,q}\}$ is the corresponding set of unknown quadrature weights.

The unknown weights are determined by enforcing the exactness conditions. Note that due to the local support of B-splines,  $B_i$ only overlaps with $2p+1$ neighboring B-splines (including itself), namely $B_j$ for $j=i-p,\ldots, i+p$. Consequently, we impose the following $2p+1$ constraints to ensure the exact integration,
\begin{align}
\sum_{q=1}^{r} w_{i,q} B_{i-p}(u_{i,q}) &= \int_0^m B_{i-p}(u) \left( B_i(u) du \right), \notag \\
&\vdots \\
\sum_{q=1}^{r} w_{i,q} B_{i+p}(u_{i,q}) &=\int_0^m B_{i+p}(u) \left( B_i(u) du \right) , \notag
\label{WQ_B}
\end{align}
where $r$ denotes the number of quadrature points located in the support of $B_i$. Our objective is to solve this local linear system and find all $w_{i,q}$. Problem (4) is well-posed when the number of unknown weights is no less than that of the exactness constraints, i.e., $r \geq 2p+1$. As a result, the linear system is generally underdetermined. It is straightforward to verify that in the uniform case, this condition is guaranteed by setting two quadrature points per knot interval. 

The underdetermined system is then solved via the quadratic optimization,
\begin{equation} \label{eq:optimization}
\begin{aligned}
\min_{\mathbf{w}} \quad & \left\| \mathbf{Z}^{-1} \mathbf{w} \right\|_{l^2}^2 \\
\text{s.t.} \quad & \mathbf{A} \mathbf{w} = \mathbf{b} \\
&\mathbf{A} = \begin{pmatrix}
B_{i-p}(u_{i,1}) & \cdots & B_{i-p}(u_{i,r}) \\
\vdots & & \vdots \\
B_{i+p}(u_{i,1}) & \cdots & B_{i+p}(u_{i,r})
\end{pmatrix}, \quad\mathbf{b}=\begin{pmatrix}
\int_0^m B_{i-p}(u) B_i(u) du \\
\vdots \\
\int_0^m B_{i+p}(u) B_i(u) du
\end{pmatrix},\\
&\mathbf{Z}= \text{diag}\left\{ B_i(u_{i,1}),\cdots,B_i(u_{i,r})\right\},
\end{aligned}
\end{equation}
where  $\left\| \cdot \right\|_{l^2}^2$  is the Euclidean norm,  $\mathbf{w}=\{w_{i,1},\ldots, w_{i,r}\}$, $\mathbf{A}\in\mathbb{R}^{(2p+1)\times r}$, $\mathbf{b}\in\mathbb{R}^{2p+1}$, and the diagonal matrix $\mathbf{Z}$ is introduced to normalize the weights. Minimizing $\left\| \mathbf{Z}^{-1} \mathbf{w} \right\|_{l^2}^2$ leads to a balanced weight distribution. The entries of $\mathbf{b}$ are precomputed \cite{lookup2015}. The constrained optimization problem is solved by the Lagrangian method, ultimately leading to
\begin{equation}
\mathbf{w}=\mathbf{Z}^2\mathbf{A}^\top(\mathbf{A}\mathbf{Z}^2\mathbf{A}^\top)^{-1}\mathbf{b},
\label{eq:w0}
\end{equation}
where the inverse of $\mathbf{A}\mathbf{Z}^2\mathbf{A}^\top = (\mathbf{A}\mathbf{Z})(\mathbf{A}\mathbf{Z})^\top$ is computed with the help of QR factorization. Letting $\mathbf{v}=\mathbf{Z}^{-1}\mathbf{w}$,  $\mathbf{Aw}=\mathbf{b}$ is equivalent to $\mathbf{(AZ)v}=\mathbf{b}$. Performing QR factorization to $(\mathbf{A}\mathbf{Z})^\top$ yields $(\mathbf{A}\mathbf{Z})^\top = \mathbf{Q}\mathbf{R}$. Substituting it into Eq.~\eqref{eq:w0}, we solve $\mathbf{w}$ with
\begin{equation}
\mathbf{w}=\mathbf{Z}^{-1}\mathbf{Q}(\mathbf{R}^\top)^{-1}\mathbf{b}.
\label{QR}
\end{equation}
Note that $\mathbf{R}^\top$ is a triangular matrix, so solving it is straightforward.

\section{Analysis-suitable unstructured T-splines (ASUTS)}
In this section, we present the basics of analysis-suitable unstructured T-splines (ASUTS). Among many other options for unstructured splines, ASUTS possesses desired properties for both design and analysis, such as polynomial partition of unity, refinability, optimal convergence, etc. On the other hand, its construction is fairly complicated, which involves mixed types of splines and a dedicated refinement scheme that passes along the essential information across refinement levels. We believe that it will serve as a non-trivial foundation for our discussion. Interested readers may refer to \cite{Wei2022AnalysissuitableUT,Casquero2020SeamlessIO} for details.

ASUTS has two building blocks: a control mesh and associated splines. The control mesh, comprising vertices (control points), edges, and faces (elements), typically takes the form of an unstructured quadrilateral (quad) mesh that may contain T-junctions. T-junctions are introduced to support adaptive refinement, but they are not directly relevant to this work and thus omitted from our subsequent discussion. An example of a control mesh of interest is shown in Fig. 1. We first introduce several terminologies.

\begin{itemize}
\item \textit{Valence} of a point: the number of edges sharing the point. 
\item \textit{Extraordinary point} (EP): an interior point with valence other than four, or a boundary point with valence exceeding three. 
\item \textit{Spoke edge}: an edge with a vertex being an EP. 
\item \textit{Irregular element}: an element containing an EP. 
\item \textit{Enriched element}: an element where extra degrees of freedom (DOF) are introduced. This includes all the irregular elements in the initial control mesh, and the child elements of every enriched element after refinement.
\item \textit{Transition point}: a point shared by both an enriched element and a non-enriched element, i.e., a point lying on the interface between enriched elements and non-enriched elements.
\item \textit{Transition element}: an element containing a transition point. It can be either enriched or non-enriched.

\item \textit{Neighborhood}: the 1-ring neighborhood of a vertex consists of all the elements sharing the vertex. The $n$th-ring neighborhood of a vertex comprises all elements sharing vertices with its ($n-1$)-ring neighborhood for $n\geq 2$. The $n$-ring neighborhood of a vertex is the union of its ($n-1$)-ring and $n$th-ring neighborhoods. The $n$-ring neighborhood of the boundary is the union of the $n$-ring neighborhoods of all the boundary vertices. The $n$th-ring neighborhood of the boundary is the difference set between the $n$-ring and ($n-1$)-ring neighborhoods.
\item \textit{Boundary element}: an element in the 2nd-ring neighborhood of the boundary. It can be enriched or non-enriched. The implication of this definition will become clear when we discuss the
assumptions on the quad mesh. Elements in the 1-ring neighborhood of the boundary are distinguished as passive boundary elements since they do not contribute to analysis.
\item \textit{Regular element}: an element that is not enriched, transition, or boundary.
\item Knot span of an edge: a non-negative parameter assigned to each edge for the parameterization purpose. The knot spans of opposite edges in an element must be identical. They determine the parametric area of an element.
\end{itemize}

We adopt the following mild assumptions on the input control mesh to simplify the discussion while maintaining a reasonably general setting.
\begin{itemize}
\item EPs are not allowed in the 1-ring neighborhood of the boundary. In other words, all the vertices in the 1-ring neighborhood of the boundary are regular.
\item All knot spans are set to be 1 except for those perpendicular to the boundary, which are 0. This setting enables the treatment of boundary using open knot vectors for bicubic splines. Consequently, all passive boundary elements have a zero parametric area, whereas the elements in the 2nd-ring neighborhood of the boundary are the effective boundary elements in analysis.
\end{itemize}

\begin{figure}[htbp]
\parbox{0.45\textwidth}{\centering\includegraphics[width=\linewidth]{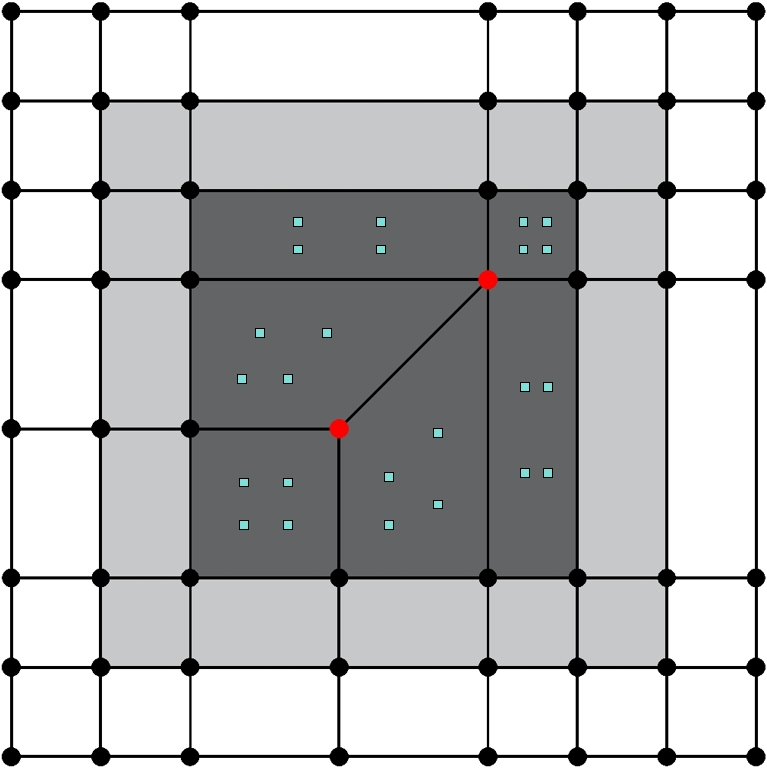}\\(a)} 
\hfill 
\parbox{0.45\textwidth}{\centering\includegraphics[width=\linewidth]{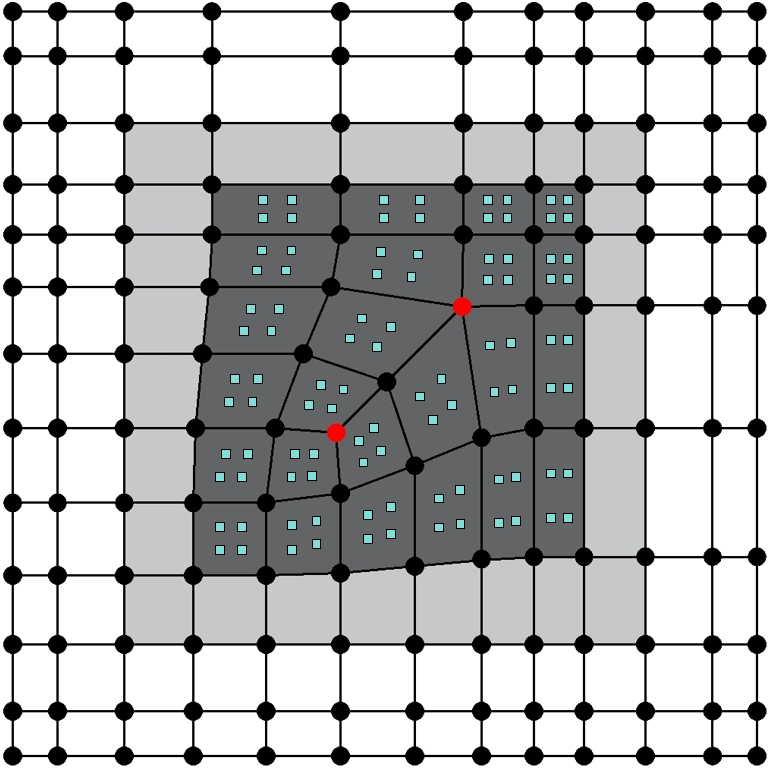}\\(b)} 
\caption{Illustration of ASUTS control meshes.  (a) An initial control mesh with two extraordinary points (marked as red dots) and four face-based control points (green squares) in each irregular/enriched face. (b) The control mesh after global refinement once according to ASUTS refinement rules. Enriched and transition faces are shaded dark and light gray, respectively.}
\label{fig:asuts_control_meshes} 
\end{figure}

The basis functions of ASUTS are defined elementwise. They are simply bicubic B-splines in regular elements. On the other hand, irregular and transition elements need special treatment to ensure important properties such as refinability, partition of unity, and optimal convergence.

In irregular elements, ASUTS employs the D-patch method \cite{REIF1997174}, where degenerate Bézier patches join smoothly around EPs without cusps. Without loss of generality, we focus on the \textit{analysis space} of ASUTS in this work, where 4 additional face-based control points are introduced in each irregular/enriched element,  whereas the splines associated with EPs are eliminated; see Fig. 1. The face-based points are associated with $C^1$ bicubic B-splines, each of whose knots is repeated twice. This enrichment strategy is essential for guaranteeing refinability and optimal convergence.

Each irregular element is first converted to a $C^0$ Bézier patch, followed by a $2\times 2$ split to increase the available degrees of freedom (DOFs) to enforce continuity constraints; see Fig. 2. Consequently, each spline $N_i$ of ASUTS is defined piecewise with respect to these sub-elements. In each sub-element, \(N_i\) is expressed as a linear combination of $4 \times 4$ bicubic Bernstein polynomials, with specific coefficients  determined by the D-patch method to enforce $C^1$ continuity across spoke edges. In Fig. 2, DOFs indicated by square markers are collapsed to identical positions at the extraordinary point, whereas those marked by triangles are made co-plane. DOFs corresponding to solid disks are further adjusted through a simple averaging of neighboring DOFs. The remaining coefficients (hollow circles) are computed using the classical de Casteljau's algorithm~\cite{Piegl2012}. 

\begin{figure}[t]
\centering
\includegraphics[width=7cm]{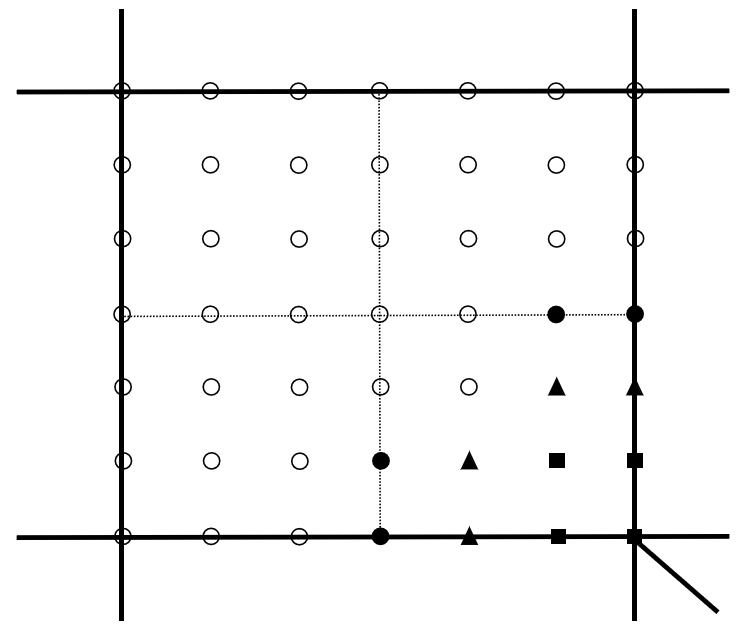}
\caption{Spline construction within an irregular element using the D-patch method \cite{REIF1997174}. Squares are collapsed to identical positions at the EP, triangles are made co-plane, hollow circles correspond to coefficients from the classical de Casteljau algorithm, and coefficients of solid disk are determined through averaging with neighboring triangles and circles.}
\label{Dpatch}
\end{figure}

A transition element involves both vertex-based DOFs (black dots in Fig. 1) and face-based DOFs. These two distinct spline types are blended together using the truncation mechanism originally proposed in \cite{WEI2018609}. This approach eliminates redundant contributions from vertex-based splines such that they form a  partition of unity when combined with face-based splines.

To achieve optimal convergence, the irregular region must be kept the same size during refinement; see Fig. 1(a, b). This constraint ensures that ASUTS  spaces are either nested or nearly nested. During refinement, vertex-based control points are computed by the classical knot insertion algorithm, whereas the face-based control points are computed using the de Casteljau algorithm. 

\section{Weighted quadrature on ASUTS}
In this section, we introduce the proposed weighted  quadrature (WQ) method on analysis-suitable unstructured T-splines (ASUTS). We develop specialized strategies to address the unique challenges caused by EPs, including the absence of a global parametric domain, the variable number of basis functions across different element types, and the ill-conditioning of the resulting underdetermined linear systems. Note that the present work is restricted to three-dimensional surfaces (2-manifolds) and bicubic splines only. Extensions to volumes and other degree cases are possible but lie outside the scope of this work.

\subsection{Quadrature points in the physical domain}
The standard WQ on B-splines requires a global parametric domain such that the geometric and/or  PDE data, i.e.,  $c(u)$  in Eq. (2), can be parameterized in a highly smooth manner. However, unstructured splines do not have such a global parametric domain. In general, every element has its own independent local parametric domain. These local domains are fundamentally incompatible around EPs, where a common parametric coordinate system cannot be consistently defined. As a result, 
$c(u)$ exhibits unavoidable discontinuities across element boundaries due to rigid rotations and translations of the local parametric domains, which violates the fundamental assumption of the standard WQ that $c(u)$ is sufficiently smooth.

To tackle this issue, we propose to develop the WQ rule directly in the physical domain, which is clearly shared by all the elements. We begin with mass matrix terms,
\begin{equation}
 \int_\Omega c(\mathbf{x}) N_i(\mathbf{x}) N_j(\mathbf{x}) d\mathbf{x},
\end{equation}
where $\Omega$ is the computational domain parameterized by ASUTS, $c(\mathbf{x})$ comes from the PDE data, and $N_i(\mathbf{x})$ and $N_j(\mathbf{x})$ are test and trial functions, respectively. Note that $N_i(\mathbf{x})$ and $N_j(\mathbf{x})$ are bivariate ASUTS basis functions. Unlike the standard WQ where weights can be computed independently for each parametric direction, in the unstructured setting, weights must be computed for each bivariate test function as a whole due to the absence of the tensor-product structure. 

In this setting, $c(\mathbf{x})$ is smooth and thus can be neglected when defining the quadrature rule. We have
\begin{equation}
\int_{\Omega}N_i(\mathbf{x})N_j(\mathbf{x})d\mathbf{x}\approx\mathbb{Q}_{i}(N_j(\mathbf{x}))=\sum_{q=1}^{r_i}w_{i,q}N_j(\mathbf{x}_{i,q}), \quad i, j = 1, \ldots, n,
\label{WQ_define}
\end{equation}
where each quadrature point $\mathbf{x}_{i,q}$ is predefined and located in the physical domain, $w_{i,q}$  is the unknown quadrature weight associated with $\mathbf{x}_{i,q}$ for test function $N_i$, $r_i$ is the number of quadrature points for $N_i$, $N_j$ is a trial function, and $n$ is the total number of spline basis functions. The specific locations of $\mathbf{x}_{i,q}$ are not important. For implementation ease in the unstructured setting, they are evenly spaced within the parametric domain of each element and mapped to the physical domain. 

To determine $w_{i,q}$, exactness conditions must be enforced for every $N_j$ (including $N_i$ itself) whose support overlaps with that of $N_i$, that is, for splines in
\begin{equation}
\mathcal{S}_i = \{ N_j: \mathrm{supp}(N_j) \cap \mathrm{supp}(N_i) \neq \varnothing \} ,
\label{region}
\end{equation}
where $\mathrm{supp}(N_i)$ is defined as a set of element indices in ASUTS. Let $n_i$ denote the dimension of $\mathcal{S}_i$, i.e., $n_i=\#\mathcal{S}_i$. We require $r_i \geq n_i$ to ensure that the problem is well-posed. We will introduce the specific strategy to guarantee this condition in Section 4.2.

\begin{remark}
Defining quadrature rules directly in the physical domain may appear counterintuitive at the first glance, but it serves effectively the role in reducing the number of quadrature points, which is the primary focus of the work.  An immediate consequence is that the proposed approach can only be applied to problems with fixed physical domains. Any change in the physical domain will need an update for the quadrature weights. Nonetheless, the approach can be highly efficient to problems where the physical domain remains fixed and repeated computations are necessary, such as the nonlinear heat transfer problem and topology optimization.
\end{remark}

In the case of stiffness matrix terms, they involve derivatives of basis functions (e.g., gradients and Laplacian). These derivatives must be computed with respect to the physical coordinates. In particular, we consider the following terms that appear in Poisson's problem and the biharmonic problem,
\begin{equation}
\begin{aligned} 
&\quad \int_{\Omega} \nabla N_i(\mathbf{x}) \cdot \nabla N_j(\mathbf{x}) d\mathbf{x}\\ &\approx \mathbb{Q}_i^{(1),x}\left(\frac{\partial N_j(\mathbf{x})}{\partial x}\right)+\mathbb{Q}_i^{(1),y}\left(\frac{\partial N_j(\mathbf{x})}{\partial y}\right)+ \mathbb{Q}_i^{(1),z}\left(\frac{\partial N_j(\mathbf{x})}{\partial z}\right)\\ 
&=\sum_{q=1}^{r_i} w_{i,q}^{(1),x}\frac{\partial N_j(\mathbf{x}_{i,q})}{\partial x}+ w_{i,q}^{(1),y}\frac{\partial N_j(\mathbf{x}_{i,q})}{\partial y}+w_{i,q}^{(1),z}\frac{\partial N_j(\mathbf{x}_{i,q})}{\partial z}, 
 \end{aligned}
 \label{diff}
\end{equation}
\begin{equation}
 \quad \int_{\Omega} \Delta N_i(\mathbf{x}) \Delta N_j(\mathbf{x}) d\mathbf{x} \approx \mathbb{Q}_i^{(2)}(\Delta N_j(\mathbf{x})) = \sum_{q=1}^{r_i} w_{i,q}^{(2)} \Delta N_j(\mathbf{x}_{i,q}) ,
\end{equation}
where quadrature points $\mathbf{x}_{i,q}$ are identical to those defined in Eq. (\ref{WQ_define}), and weights $w_{i,q}^{(1),x}$, $w_{i,q}^{(1),y}$, $w_{i,q}^{(1),z}$ and $w_{i,q}^{(2)}$ are the unknowns to be determined.

In the case of 3D surfaces (more precisely, 2-manifold), the gradient operator $\nabla(\cdot)$ and Laplacian operator $\Delta(\cdot)$ in Eqs. (11, 12) are in fact the surface gradient operator and the Laplace–Beltrami operator \cite{Bartezzaghi2015IsogeometricAO}, respectively. Let $x(\bm{\theta})$ be a parametric surface, where $\bm{\theta}=(\theta^1,\theta^2)$ denotes the local parametric coordinates of each element, and $\mathbf{x}\in \mathbb{R}^3$ is the position vector in the 3D Euclidean space. The tangent vectors to the surface are given by
\begin{equation}
\mathbf{a}_{\alpha}=\frac{\partial\mathbf{x}}{\partial\theta^\alpha}, \quad\alpha \in \{1,2\}.
\end{equation}
The covariant metric tensor is given by
\begin{equation}
a_{\alpha\beta}=\mathbf{a}_{\alpha}\cdot \mathbf{a}_{\beta},\quad \alpha,\beta \in \{1,2\}.
\label{eq:covariant}
\end{equation}
The corresponding contravariant metric tensor is obtained by inverting the covariant metric tensor,
\begin{equation}
[a^{\alpha\beta}]=[a_{\alpha\beta}]^{-1}.
\label{eq:contravariant}
\end{equation}
The surface gradient operator and the Laplace–Beltrami operator are defined as
\begin{align}
&\nabla(\cdot)=\frac{\partial x_i}{\partial \theta^{\alpha}}\,a^{\alpha\beta}\,\frac{\partial (\cdot)}{\partial \theta^{\beta}},\\
&\Delta(\cdot)=\frac{1}{J}\frac{\partial}{\partial \theta^{\alpha}}(J\,a^{\alpha\beta}\frac{\partial (\cdot)}{\partial \theta^{\beta}}),
\end{align}
where the Einstein summation convention is employed for indices $\alpha$ and $\beta$ and $J$ denotes the determinant of the covariant metric tensor.

\subsection{Determination for the number of quadrature points}
Unlike B-splines, the number of basis functions varies across different types of elements in ASUTS due to the presence of EPs and the mixed types of splines; see Table~1. This complexity introduces a challenge when we aim to find the smallest number of quadrature points while ensuring the well-posedness condition for every test function. Our objective is to develop a general strategy applicable to all the quad meshes that satisfy the mild assumptions in Section 3.

\begin{table}[!htbp]
\centering
\caption{The number of  basis functions per element type in ASUTS across global  refinement levels, where $\mu$ denotes the valence of the extraordinary point (EP) in an irregular element and ``NA'' means that the corresponding element type does not exist at the refinement level of interest.}
\label{tab:basis_num_merged}
\begin{tabular}{p{4.3cm}p{2.5cm}p{2.5cm}p{2.8cm}}
\hline\noalign{\smallskip}
Element type & Initial mesh & 1  refinement & $\geq$2  refinements \\
\noalign{\smallskip}\svhline\noalign{\smallskip}
Regular & 16 & 16 & 16 \\
Boundary & 16 & 16 & 16 \\
Enriched regular & NA & NA & 16 \\
Non-enriched transition & 16/17 & 16/17 & 16/17 \\
Enriched transition & $3\mu+13$ & 20/21 & 20/21 \\
Irregular & $3\mu+13$ & $3\mu+8$ & $3\mu+8$ \\
\noalign{\smallskip}\hline\noalign{\smallskip}
\end{tabular}
\end{table}

\begin{figure}[htbp]
    \centering
    \begin{subfigure}{0.32\textwidth}
        \centering
        \includegraphics[width=\linewidth]{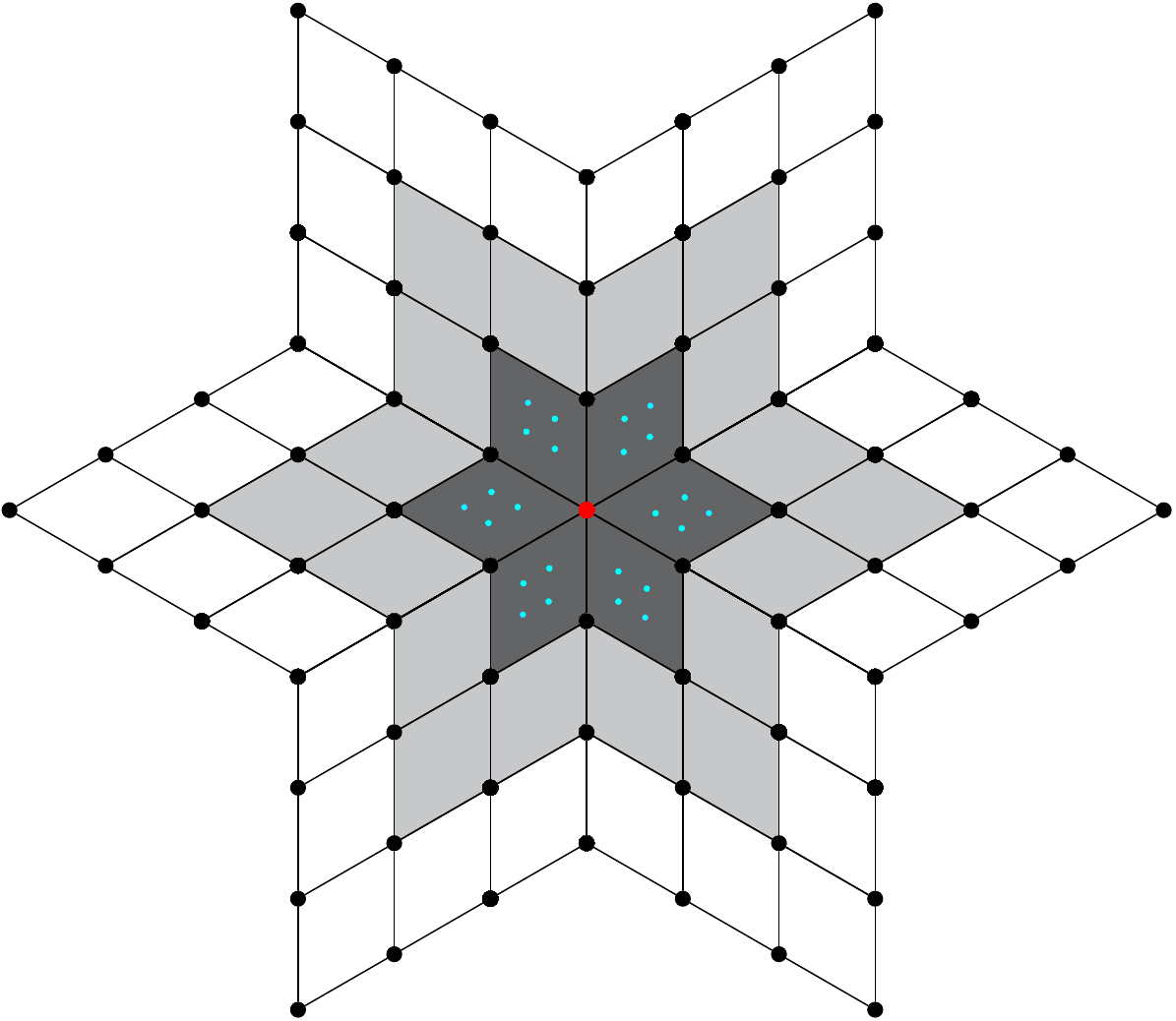}
        \caption{Initial mesh} 
    \end{subfigure}
    \hfill
    \begin{subfigure}{0.32\textwidth}
        \centering
        \includegraphics[width=\linewidth]{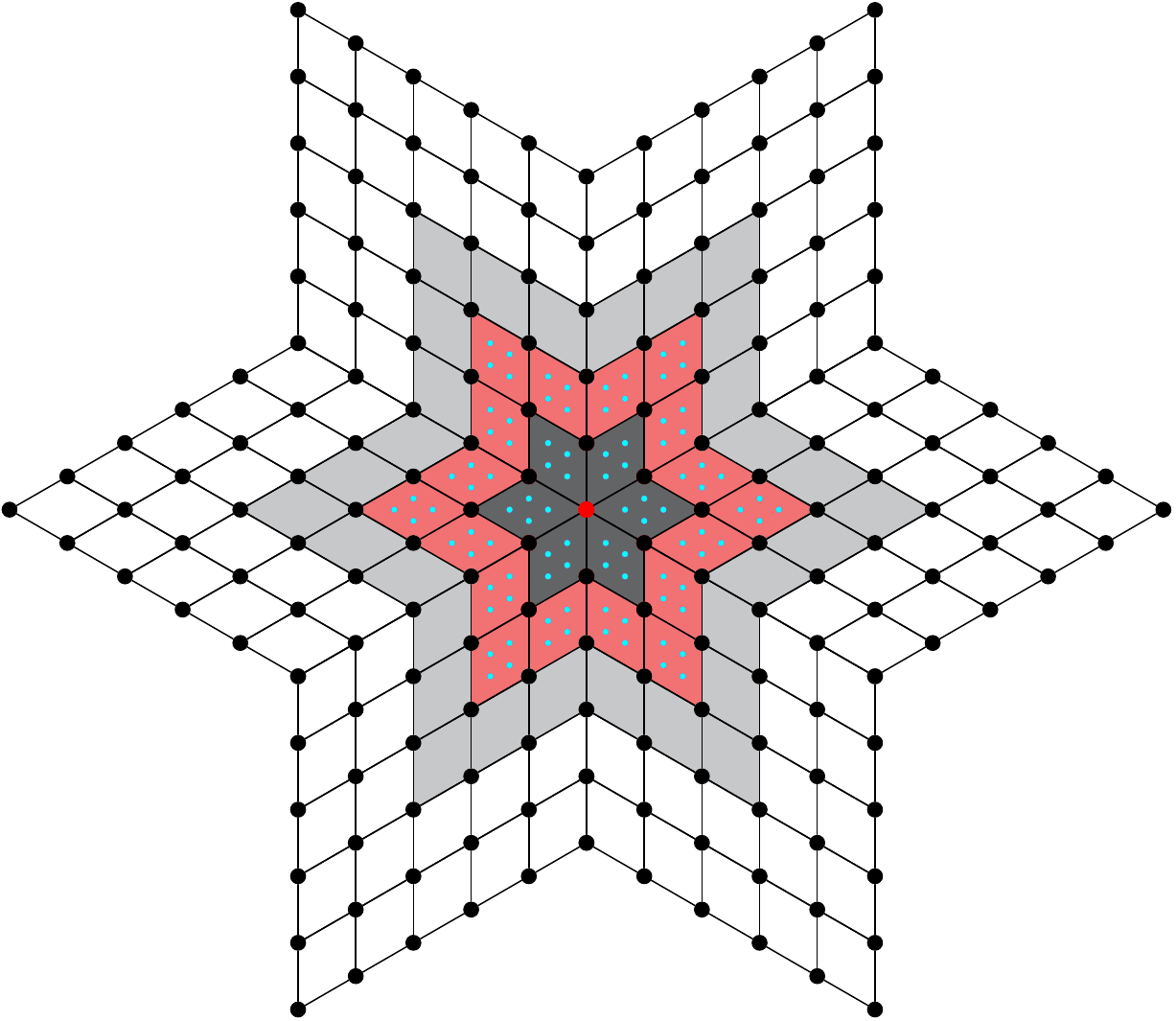}
        \caption{Global refinement once}
    \end{subfigure}
    \hfill
    \begin{subfigure}{0.32\textwidth}
        \centering
        \includegraphics[width=\linewidth]{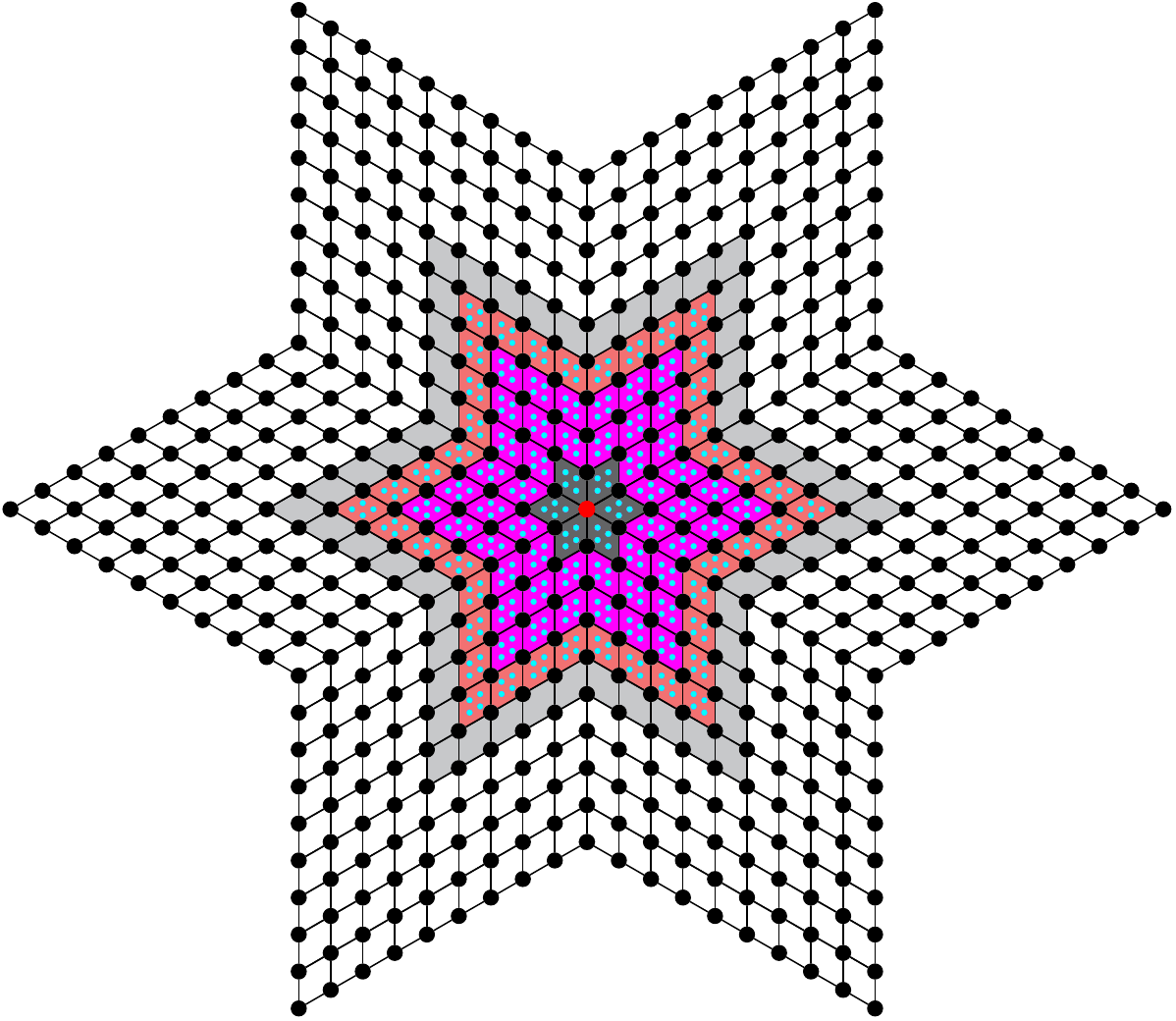}
        \caption{Global refinements twice}
    \end{subfigure}
    \caption{Element types at different refinement levels. Distinct colors represent different element types: dark gray for irregular elements, light pink for enriched transition elements, magenta for enriched regular elements, light gray for transition elements, and white for regular elements. The red dot is an EP, whereas green dots are face-based points.}
    \label{fig:local_refine_element}
\end{figure}

We first comment on the numbers in Table~1. In the initial mesh, irregular elements are also enriched elements and transition elements. A non-enriched transition element may have either 16 or 17 basis functions, depending on whether it shares a vertex (16) or an edge (17) with an enriched element. After one global refinement, some enriched elements are no longer irregular because they do not involve EPs; see the 2nd-ring neighborhood of the EP in Fig. \ref{fig:local_refine_element}(b).  Instead, they are transition elements and such an element may have either 20 or 21 basis functions, depending on whether it shares an edge (20) or a vertex (21) with a non-enriched element. After two (or more) global refinements, some enriched elements are neither irregular nor transition elements, and thus called \textit{enriched regular} elements. For example, the 2nd- and 3rd-ring elements of an EP fall into this category; see Fig. \ref{fig:local_refine_element}(c). 

We propose a ``min-max" strategy to determine the number of quadrature points element-wise such that the well-posedness is guaranteed for every test function. First, for a certain test function $N_i$ and for every element in its support, $\forall k\in\mathrm{supp}(N_i)$, we find the minimum number of quadrature points $r_{i,k}$ for element $k$ such that overall it guarantees $\sum_{k} r_{i,k} \geq n_i$, where recall that $n_i = \#\mathcal{S}_i$; see Eq. (10). $r_{i,k}$ is determined by the element type; see Table~2 for the final result.

\begin{table}[htbp]
\centering
\caption{The number of weighted quadrature points per element type, with the standard Gaussian quadrature (GQ) shown as a reference, where ``QP" stands for ``Quadrature Points".}
\label{tab:quadrature_distribution}
\begin{tabular}{p{4.5cm}p{2.8cm}p{4.2cm}}
\toprule
Element type & \#QP (WQ) & \#QP (GQ) \\
\midrule
Regular                        & $2\times 2$       & $4\times 4$ \\
Boundary                       & $4\times 4$       & $4\times 4$ \\
Irregular                      & $4\times 4$       & $4\times 4\times 4$ \\
Enriched transition            & $4\times 4$  & $4\times 4$ \\
Non-enriched transition        & $3\times 3$ / $4\times 4$ & $4\times 4$ \\
Enriched regular               & $3\times 3$ / $4\times 4$ & $4\times 4$ \\
\bottomrule
\end{tabular}
\\[6pt]
\raggedright
\footnotesize\textit{Note:} The particular choice of ``$3\times 3$ / $4\times 4$''  depends on the valence and number of involved EPs. Specifically, $4\times 4$ quadrature points are required in the case of high-valence EPs (valence~$\geq 7$) and clustered EPs, whereas $3\times 3$ quadrature points are sufficient in the other cases.
\end{table}

        
        

In Table~2, the numbers corresponding to regular and boundary elements follow the standard WQ to stay compatible. Next, the number for the other element types is determined proportionally to the number $(m_k)$ of basis functions defined on the element of interest. Specifically, we have
\begin{equation}
r_{i,k} = \frac{m_k}{\sum\limits_{l \in \text{supp}(N_i) \setminus \mathcal{R}} m_l} \left( n_i - e_i^r \times 4 - e_i^b \times 16 \right),
\label{eq:qp_allocation}
\end{equation}
where $\mathcal{R}$ is the index set of regular elements and boundary elements, and $e_i^r$ and $e_i^b$ are the number of regular and boundary elements within the support of $N_i$, respectively. This setting is motivated by the observation that more quadrature points are needed if the element of interest has more supported functions.

After traversing every test function, every element will have one or more required numbers of quadrature points. Among them, we pick the maximum number, i.e., for element $k$, the number is given by $r_k = \max_i r_{i,k}$.

By performing this ``min-max" strategy for various inputs, we find that the numbers reported in Table~2 always satisfy the well-posed condition. However, if these numbers violate the condition in some cases, we can simply return to the general ``min-max" strategy to adjust the required numbers. By construction, this strategy naturally guarantees the well-posed condition for ASUTS. It may also be applicable to other types of unstructured splines and other degree cases.

Compared with GQ, the proposed WQ scheme achieves a significant reduction in the total number of quadrature points.  Note that an irregular element undergoes a $2\times 2$ split, where every sub-element needs $4\times 4$ Gaussian points, so in total an irregular element needs 64 Gaussian points. The reduction in the number of quadrature points is particularly significant when regular elements dominate the control mesh. This can be leveraged with sophisticated quad meshing tools that feature a small number of EPs, such as \cite{Wei2025}.

\begin{remark}
A general input control mesh may have clustered EPs in a local region. As a result, additional face-based DOFs are introduced, leading to the increase in the number of function overlaps (i.e., $n_i$). Since the ratio $m_k/ \sum m_l$ in Eq. (\ref{eq:qp_allocation}) remains approximately constant, the increase in $n_i$ directly yields a larger required number of quadrature points $r_{i,k}$. Consequently, elements near clustered EPs demand more quadrature points; see the ``$4 \times 4$" option in Table~2 for non-enriched transition elements and enriched regular elements.
\end{remark}

\subsection{Treatment of ill-conditioned underdetermined systems}
The well-posedness condition generally results in an underdetermined local linear system for each test function. In particular, splines surrounding extraordinary points frequently lead to numerically ill-conditioned systems, especially when derivative operators are involved. Consequently, the computed weights exhibit extremely low accuracy if the system is solved using the standard approach described in Section 2.

We address this issue by improving the conditioning of the coefficient matrix. The treatment is identical for different kinds of matrices, so we will use the stiffness matrix of Poisson's problem to explain the idea. Similarly to Eq. (5), the coefficient matrix of the underdetermined linear system is

\begin{equation}
\mathbf{A} = 
\begin{bmatrix} 
\dfrac{\partial N_1}{\partial x}(\mathbf{x}_{i,1}) & \cdots & \dfrac{\partial N_1}{\partial x}(\mathbf{x}_{i,r_i}) \\ 
\vdots & & \vdots \\ 
\dfrac{\partial N_{n_i}}{\partial x}(\mathbf{x}_{i,1}) & \cdots & \dfrac{\partial N_{n_i}}{\partial x}(\mathbf{x}_{i,r_i})
\end{bmatrix},
\label{system}
\end{equation}

We first perform singular value decomposition (SVD) on $\mathbf{A}$:
\begin{equation}
\mathbf{A}=\mathbf{U\Sigma V}^\top,
\end{equation}
where $\mathbf{U} \in \mathbb{R}^{n_i\times n_i}$, $\mathbf{V} \in \mathbb{R}^{r_i\times r_i}$ are orthogonal matrices, and $\mathbf{\Sigma} \in \mathbb{R}^{n_i\times r_i}$ is a diagonal matrix with positive entries $\sigma_k$ referred to as singular values. Arranging $\sigma_k$ in descending order, we truncate $\mathbf{A}$ by discarding extremely small singular values. Such singular values correspond to unstructured splines around extraordinary points, typically on the order of $10^{-16}$, while the other singular values are greater than $10^{-6}$. Equivalently, we retain the first $t$ singular values, and $\mathbf{A}$ is approximated as
\begin{equation}
\mathbf{A}\approx \mathbf{U}_t\mathbf{\Sigma}_t\mathbf{V}_t^\top,
\end{equation}
where $\mathbf{U}_t\in\mathbb{R}^{n_i\times t}$ contains the first $t$ columns of $\mathbf{U}$, $\mathbf{V}_t\in\mathbb{R}^{r_i\times t}$ contains the first $t$ columns of $\mathbf{V}$, and $\mathbf{\Sigma}_t\in\mathbb{R}^{t\times t}$ is the diagonal matrix of the first $t$ singular values from $\mathbf{\Sigma}$. After performing this truncated SVD \cite{FASSINO2022112746}, the weight vector $\mathbf{w}$ is computed as
\begin{equation}
\mathbf{w}=\mathbf{Z}^2\mathbf{V}_t(\mathbf{V}_t^{\top}\mathbf{Z}^{2}\mathbf{V}_t)^{-1}\mathbf{\Sigma}_t^{-1}\mathbf{U}^{\top}_t\mathbf{b}.
\end{equation}

Note that after discarding extremely small singular values, inverting $\mathbf{\Sigma}$ no longer presents numerical difficulties. Moreover, $\mathbf{V}_t^\top\mathbf{Z}^2\mathbf{V}_t$ inherits the favorable conditioning properties of $\mathbf{Z}$, enabling stable and accurate computation of quadrature weights.

\begin{remark}
While SVD is computationally costly, it is not a bottleneck in our framework for two reasons. First, quadrature weights for fixed domains are computed only once and can be reused in nonlinear problems. Second, truncated SVD only applies locally to EP-near test functions, and the number of EPs are far fewer than that of total DOFs, resulting in negligible overhead.
\end{remark}

\section{Application in the nonlinear  heat transfer problem}
In this section, we employ the proposed WQ on ASUTS to solve the nonlinear heat transfer problem, which requires repeated matrix assembly. We begin by introducing the weak form of the problem, followed by the spatial discretization using ASUTS and the temporal discretization via the generalized trapezoidal rule \cite{CORNEJOFUENTES2023116157}. The linearized form is then derived for 3D surfaces, on which we apply the proposed WQ.

\subsection{Weak form and spatial discretization}
We consider the nonlinear  heat transfer problem defined on a parametric surface, where the physical domain \(\Omega\) is a 2-manifold with boundary \(\partial\Omega\). Let \(\mathbf{x}\) denote the spatial point on \(\Omega\), and \(T(\mathbf{x},t)\) be the temperature field defined in the space-time domain \(\Omega \times [0,T_f]\). The strong form of the problem is stated as follows:
\begin{align}
&\frac{\partial T(\mathbf{x},t)}{\partial t} = \nabla \cdot \left[ k(T(\mathbf{x},t), \mathbf{x}) \nabla T(\mathbf{x},t) \right] + f(\mathbf{x},t) &&\quad  \text{in } \Omega \times (0,T_f], \notag\\
&T(\mathbf{x},t) = g(\mathbf{x},t) &&\quad  \text{on } \partial\Omega_D \times (0,T_f],  \\
&T(\mathbf{x},0) = T_0(\mathbf{x}) &&\quad \forall \mathbf{x} \in \bar{\Omega}, \notag
\end{align}
where \(g(\mathbf{x},t)\) is the prescribed temperature on the Dirichlet boundary,  \(T_0(\mathbf{x})\) is the initial temperature distribution over the closed domain $\bar{\Omega}$ and \(k(T(\mathbf{x},t),\mathbf{x})\) is the thermal conductivity, which constitutes the source of nonlinearity in the problem.

We introduce $\mathcal{S}$ and $\mathcal{V}$ as the trial and test function spaces, respectively, where \(\mathcal{S} = \{ T \in H^1(\Omega) \mid T = g \text{ on } \partial\Omega \}\) and \(\mathcal{V} = \{ w \in H^1(\Omega) \mid w = 0 \text{ on } \partial\Omega \}\). The variational formulation of the problem is to find $T\in \mathcal{S}$ such that
\begin{equation}
\label{eq:heat_residual}
R(T, w) = \int_\Omega \dot{T} w \, d\Omega + \int_\Omega k(T) \nabla T \cdot \nabla w \, d\Omega - \int_\Omega f w \, d\Omega = 0 
\end{equation}
for any  $w\in\mathcal{V}$,  where \(\dot{T} = \partial T/\partial t\) denotes the time derivative of temperature.

 We approximate the temperature field \(T(\mathbf{x},t)\) and the test function \(w(\mathbf{x})\) by the finite-dimensional subspace spanned by the ASUTS basis functions \(N_A(\mathbf{x})\). Specifically, we find \(T_h \in \mathcal{S}_h \subset \mathcal{S}\) such that  \(\forall w_h \in \mathcal{V}_h \subset \mathcal{V}\),
\begin{equation}
R(T^h, w^h) =  0, \label{eq:heat_residual_h}
\end{equation}
where $T^h$ and $w^h$ are defined as
\begin{align}
&T_h(\mathbf{x},t) = \sum_{A=1}^{n_{\text{dof}}} T_A(t) N_A(\mathbf{x}),\\
&w_h(\mathbf{x},t) = \sum_{A=1}^{n_{\text{dof}}} w_A(t) N_A(\mathbf{x}).
\end{align}
\(N_A(\mathbf{x})\) is an ASUTS basis function, \(n_{\text{dof}}\) is the total number of basis functions, and \(T_A(t)\) and \(w_A(t)\) are control variables at time $t$. We define the residual vector as 
\begin{equation}
\mathbf{R}=\{R_A\},
\end{equation}
where $R_A$ is defined as 
\begin{equation}
R_A(T, N_A) = \int_\Omega \dot{T} N_A \, d\Omega + \int_\Omega k(T) \nabla T \cdot \nabla N_A \, d\Omega - \int_\Omega f N_A \, d\Omega = 0. \label{eq:heat_residual_component}
\end{equation}

\subsection{Temporal discretization and Newton linearization}
 Now we discretize the time interval as $\{t_n\}_{n=0}^{n_T}$. The control variables at $t_n$ are denoted as $\mathbf{T}_n = \left\{ T_n(t_n) \right\}_{n=1}^{n_T}$ and $\dot{\mathbf{T}}_n = \left\{ \frac{\partial T_n(t_n)}{\partial t} \right\}_{n=1}^{n_T}$. 
The discretized problem can be stated as follows: given  $\mathbf{T}_{n},\dot{\mathbf{T}}_n\; \text{and}\; \Delta t_n=t_{n+1}-t_n$, find $\mathbf{T}_{n+1}$ and $\dot{\mathbf{T}}_{n+1}$ such that 
\begin{align}
&\mathbf{R}(\mathbf{T}_{n+1},\dot{\mathbf{T}}_{n+1})=0,\label{eq:system}\\
&\mathbf{T}_{n+1} = \mathbf{T}_n + \Delta t \left[ (1-\theta) \dot{\mathbf{T}}_n + \theta \dot{\mathbf{T}}_{n+1} \right].\notag
\end{align}
Substituting the discretized problem into Eq. (\ref{eq:heat_residual_component}), the discrete residual at \(t_{n+1}\) becomes:
\begin{equation}
R_A(T_{n+1}) = \int_\Omega \dot{T}_{n+1} N_A \, d\Omega + \int_\Omega k(T_{n+1}) \nabla T_{n+1} \cdot \nabla N_A \, d\Omega - \int_\Omega f_{n+1} N_A \, d\Omega = 0, \label{eq:heat_residual_discrete}
\end{equation}
where \(f_{n+1} = f(\mathbf{x}, t_{n+1})\) is the heat source at step \(n+1\).

The nonlinear system of Eq. (\ref{eq:system}) is solved by the Newton method, which leads to a two-stage predictor–corrector algorithm \cite{CORNEJOFUENTES2023116157}.  During the predictor stage, the previous results $\mathbf{T}_n$ and $\dot{\mathbf{T}}_n$
 are used to approximate $\mathbf{T}_{n+1}$. 
During the corrector stage, we  linearize $\mathbf{R}(\mathbf{T}^{(i)}_{n+1},\dot{\mathbf{T}}^{(i)}_{n+1})$ to approximate $\mathbf{T}_{n+1}$ until convergence
is achieved. For iteration $i$, the linearized system is defined as
\begin{equation}
-\mathbf{K}_{n+1}^{(i)}\Delta \mathbf{T}_{n+1}^{(i)}=\mathbf{R}^{(i)}_{n+1},
\label{eq:newton}
\end{equation}
where the tangent matrix is defined as
\begin{align}
\mathbf{K}_{n+1} &= \frac{\partial \mathbf{R}(\mathbf{T}_{n+1}, \dot{\mathbf{T}}_{n+1})}{\partial \dot{\mathbf{T}}_{n+1}} \frac{\partial \dot{\mathbf{T}}_{n+1}}{\partial \mathbf{T}_{n+1}} + \frac{\partial \mathbf{R}(\mathbf{T}_{n+1}, \dot{\mathbf{T}}_{n+1})}{\partial \mathbf{T}_{n+1}} \notag\\
&= \frac{1}{\theta \Delta t} \frac{\partial \mathbf{R}(\mathbf{T}_{n+1}, \dot{\mathbf{T}}_{n+1})}{\partial \dot{\mathbf{T}}_{n+1}} + \frac{\partial \mathbf{R}(\mathbf{T}_{n+1}, \dot{\mathbf{T}}_{n+1})}{\partial \mathbf{T}_{n+1}}. \notag\\
\end{align}
Specifically, each entry of the tangent matrix is given by
\begin{align}
K_{AB}(T_{n+1}) &=  \frac{1}{\theta \Delta t} \int_\Omega N_A N_B \, d\Omega 
+ \int_\Omega \frac{k(T_{n+1})}{T_{n+1}} N_A \nabla T_{n+1} \cdot \nabla N_B \, d\Omega \notag\\&+ \int_\Omega k(T_{n+1}) \nabla N_A \cdot \nabla N_B \, d\Omega.
\end{align}
Note that the exact tangent matrix $\mathbf{K}_{n+1}$ is non-symmetric due to the presence of the term $N_A \nabla T_{n+1} \cdot \nabla N_B$. To enhance efficiency, Hughes et al. \cite{WINGET1985711} proposed a symmetric form of the tangent matrix by neglecting the non-symmetric term, 
\begin{equation}
S_{AB}(T_{n+1}) \approx K_{AB}(T_{n+1}) = \frac{1}{\theta \Delta t} \int_\Omega N_A N_B \, d\Omega+ \int_\Omega k(T_{n+1}) \nabla N_A \cdot \nabla N_B  \,d\Omega,
\label{eq:S}
\end{equation}
which still guarantees the convergence of Newton iterations. Applying  WQ to Eq. (\ref{eq:S}) yields the following discretized form:
\begin{align}
S_{AB} &\approx \frac{1}{\theta \Delta t} \sum^{r_A}_{q=1}w_{A,q}N_B(\mathbf{x}_{A,q})\notag\\ 
&+k(T_{n+1}(\mathbf{x}_{A,q}))\left[w_{A,q}^{(1),x}\frac{\partial N_B(\mathbf{x}_{A,q})}{\partial x}+w_{A,q}^{(1),y}\frac{\partial N_B(\mathbf{x}_{A,q})}{\partial y}+w_{A,q}^{(1),z}\frac{\partial N_B(\mathbf{x}_{A,q})}{\partial z}\right],
\end{align}
where $k(T_{n+1}(\mathbf{x}_{A,q}))$ denotes the temperature-dependent thermal conductivity, quadrature points $\mathbf{x}_{A,q}$ are determined according to Table~2, and weights $w_{A,q}, w_{A,q}^{(1),x}, w_{A,q}^{(1),y}$, and $w_{A,q}^{(1),z}$ are obtained by solving the linear systems described in Section 4.3.

\section{Numerical examples}
This section presents a variety of numerical tests to demonstrate the accuracy and efficiency of the proposed WQ method. We start with two  benchmark problems in 2D, the linear Poisson’s problem and the fourth-order biharmonic problem. Next, we study the nonlinear heat transfer problem defined on 3D complex surfaces to show the feasibility of the proposed method.

\subsection{Benchmark problems in 2D
}
\begin{figure}[htbp]
\centering
\includegraphics[width=6cm]{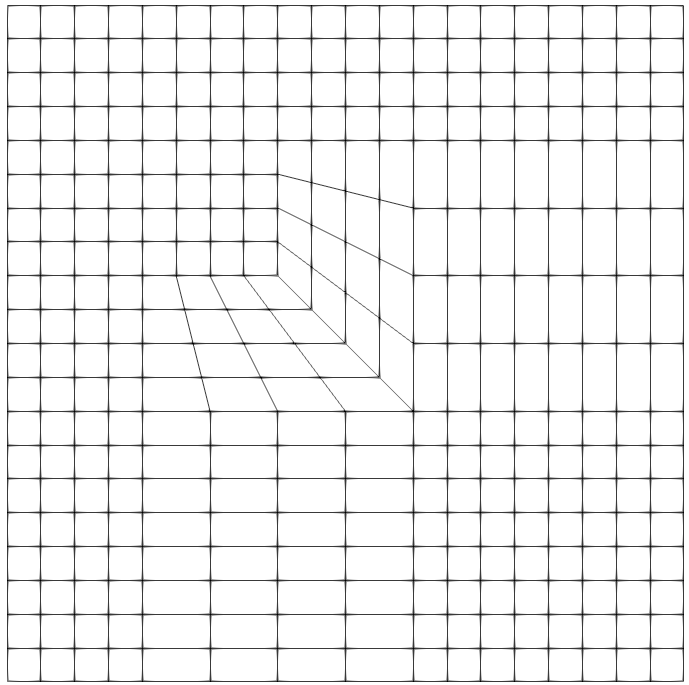}
\caption{The initial control mesh of a unit square \(\Omega = [0, 1]^2\) with two EPs.}
\label{cn}
\end{figure}

\begin{figure}[htbp]
\parbox{0.45\textwidth}{\centering\includegraphics[width=\linewidth]{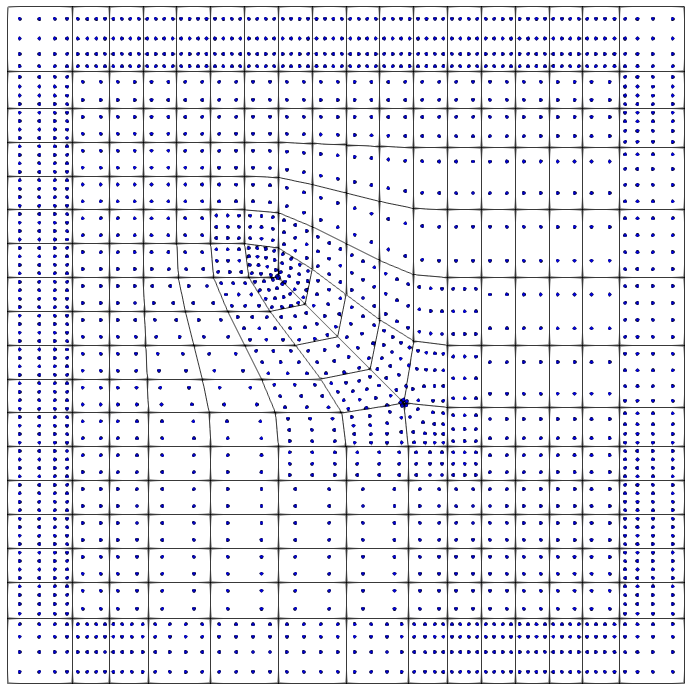}\\(a) Weighted quadrature (WQ)} 
\hfill 
\parbox{0.45\textwidth}{\centering\includegraphics[width=\linewidth]{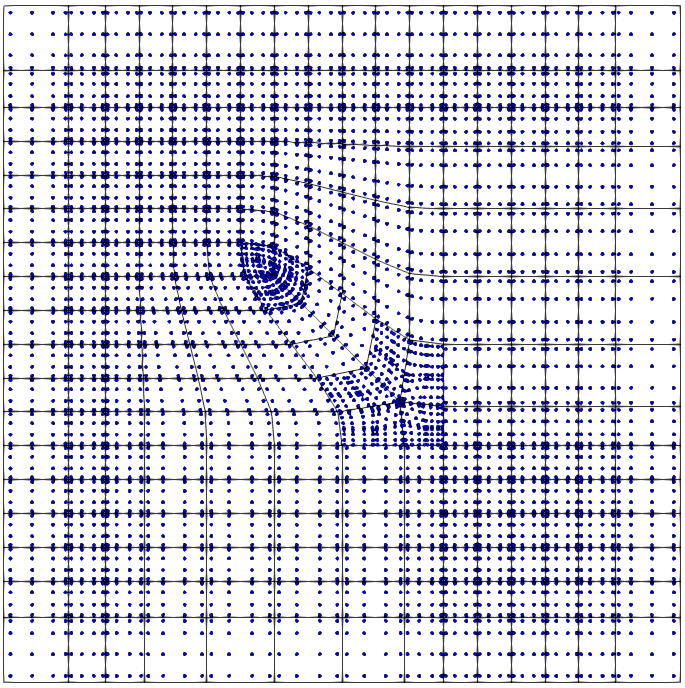}\\(b) Gaussian quadrature (GQ)} 
\caption{Locations of quadrature points in WQ and GQ, which are shown on top of B\'ezier elements.}
\label{fig:wq_gq_quad_points} 
\end{figure}

We consider a unit square domain \(\Omega = [0, 1]^2\) . It is represented using ASUTS, where two EPs are intentionally introduced in the initial control mesh; see Fig. \ref{cn}. The locations of quadrature points using WQ and Gaussian quadrature (GQ) are shown in Fig. \ref{fig:wq_gq_quad_points}, where we observe that WQ uses far fewer quadrature points than GQ. Detailed numbers are reported in Table~\ref{tab:quad_count} at different refinement levels. We observe that as the refinement level increases, regular elements become more dominant in the refined mesh, and thus fewer quadrature points are needed in WQ. In other words, the reduction ratio, computed as the relative difference between the quadrature point numbers of WQ and GQ, increases with the refinement level.

\begin{table}[htbp]
\centering
\caption{The number of quadrature points and the reduction ratio of WQ against GQ across refinement levels.}
\label{tab:quad_count}
\begin{tabular}{p{3cm}p{2.5cm}p{2.5cm}p{3.5cm}}
\hline\noalign{\smallskip}
Refinement level & WQ & GQ & Reduction\\
\noalign{\smallskip}\svhline\noalign{\smallskip}
0 (Initial) & 1944 & 4,416 &  56.0\% \\
1 & 5,936 & 16,512 & 64.0\% \\
2 & 20,240 & 64,896 & 68.8\% \\
3 & 74,744 & 258,432 & 71.0\% \\
\noalign{\smallskip}\hline\noalign{\smallskip}
\end{tabular}
\\[6pt]
\raggedright
\footnotesize\textit{Note:} The reduction ratio is computed by $\frac{\text{GQ}-\text{WQ}}{\text{GQ}}\times 100\%$.  \par

\end{table}

\begin{figure}[htbp]
\centering
\includegraphics[width=9cm]{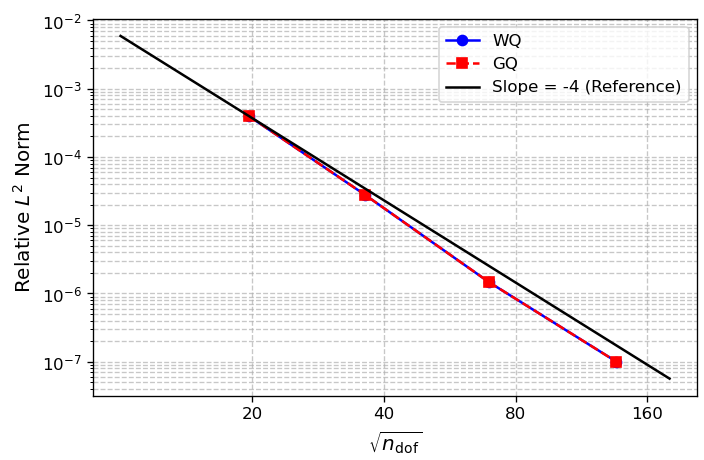}
\caption{Convergence tests of Poisson's problem on a unit square, where both WQ and GQ are used to obtain the stiffness matrix.}
\label{Poisson}
\end{figure}

We study the following Poisson's problem,
\begin{align}
-\Delta u &= f  &\text{in } &\Omega,  \notag\\
u &= 0  &\text{on } &\partial\Omega. 
\end{align}
We perform convergence tests for Poisson's problem with a manufactured solution \(u = \sin(\pi x) \sin(\pi y)\).  Three consecutive meshes are obtained by global refinement of ASUTS. The convergence curves using WQ and GQ are shown in Fig. \ref{Poisson}, where the two curves overlap one another. It means that both methods have the same level of accuracy. Indeed, they achieve the optimal convergence rate (i.e., $4$ in the $L^2$-norm). This verifies that WQ does not compromise the approximation property despite that it uses much fewer quadrature points than GQ.

We further investigate the accuracy of WQ in detail. Using GQ as the reference, entry-wise error of the stiffness matrix is computed, where the maximum absolute error is in the order of \(10^{-12}\) and maximum relative error is in the order of \(10^{-8}\) across all refinement levels. This result confirms the high level of accuracy attained by WQ.

\begin{figure}[htbp]
\centering
\includegraphics[width=9cm]{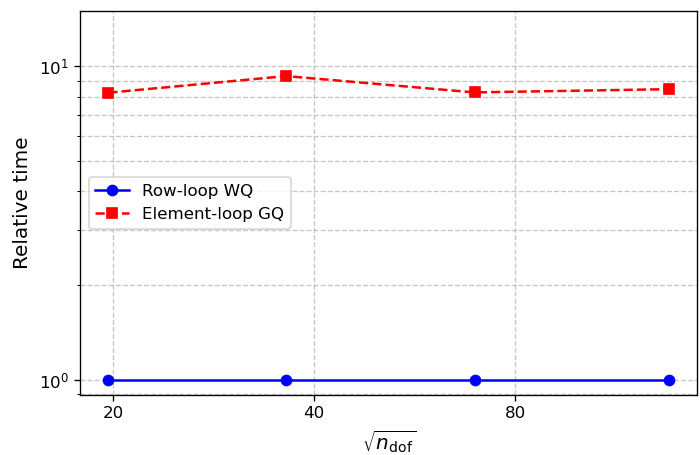}
\caption{Assembly runtime of the stiffness matrix using WQ and GQ for Poisson's problem, where the time is normalized with respect to that of WQ.}
\label{Poisson_time}
\end{figure}

Regarding efficiency, we compare the runtime of forming and assembling the stiffness matrix using WQ and GQ without sum factorization. The matrix is formed row by row using WQ, whereas it is formed element by element with GQ. To focus on the assembly efficiency, weight computation and basis function evaluation are precomputed and excluded from the reported assembly time.  Fig. \ref{Poisson_time} shows normalized runtimes with respect to that of WQ for easy comparison. We observe that GQ requires approximately 9 times of the runtime across all refinement levels.

\begin{figure}[htbp]
\centering
\includegraphics[width=9cm]{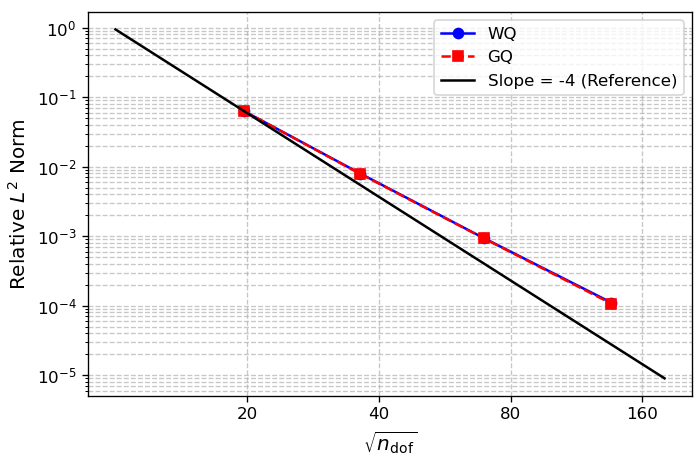}
\caption{Convergence tests of biharmonic problem on a unit square, where both WQ and GQ are used to obtain the stiffness matrix.}
\label{biharmonic}
\end{figure}

Now we consider a high-order problem, i.e., the biharmonic problem as follows,
\begin{align}
\Delta^2 u &= g  & \text{in } &\Omega,  \notag\\
u &= 0  &\text{on } &\partial \Omega,  \notag\\
\nabla u \cdot \hat{\mathbf{n}} &= 0  &\text{on } &\partial \Omega, 
\end{align}
where $\hat{\mathbf{n}}$ is the outward unit normal to the boundary. We perform convergence tests with a manufactured solution 
$v = (1 - \cos(2\pi x))$ $(1 - \cos(2\pi y))$.  Fig.  \ref{biharmonic} shows the convergence curves using WQ and GQ, which overlap with one another. This result confirms that the proposed method can also handle high-order problems without loss of accuracy. 

Regarding the entry-wise error of the stiffness matrix, the maximum absolute error is in the order of \(10^{-10}\) and the maximum relative error is in the order of \(10^{-9}\) across all refinement levels.

\begin{figure}[htbp]
\centering
\includegraphics[width=9cm]{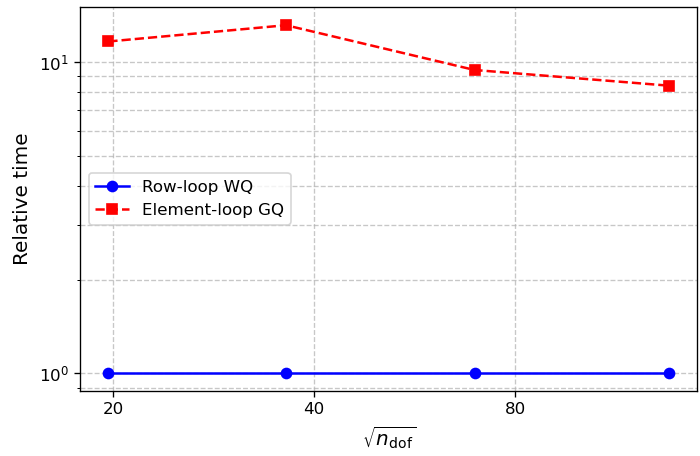}
\caption{Assembly runtime of the stiffness matrix using WQ and GQ for biharmonic problem, where the time is normalized with respect to that of WQ.}
\label{biharmonic_time}
\end{figure}

Finally, the runtime of forming and assembling the stiffness matrix using WQ and GQ is reported in Fig. \ref{biharmonic_time}. We find that WQ achieves at least 9 times speedup over GQ across all refinement levels.

\subsection{Nonlinear heat transfer problem}
Next, we study the proposed method in a more practical setting, i.e., nonlinear problems defined on 3D complex surfaces, where repeated matrix assembly usually account for 50\%$\sim$60\% of the overall simulation time. Therefore, improving the quadrature efficiency plays a key role in such problems.  For all the tests in this section,  the time integration parameter $\theta $ is set to $1$; see Eq. (30). The maximum number of Newton iterations per time step is set to $ 7$. 

\begin{figure}[htbp]
    \centering
    \includegraphics[width=0.6\textwidth]{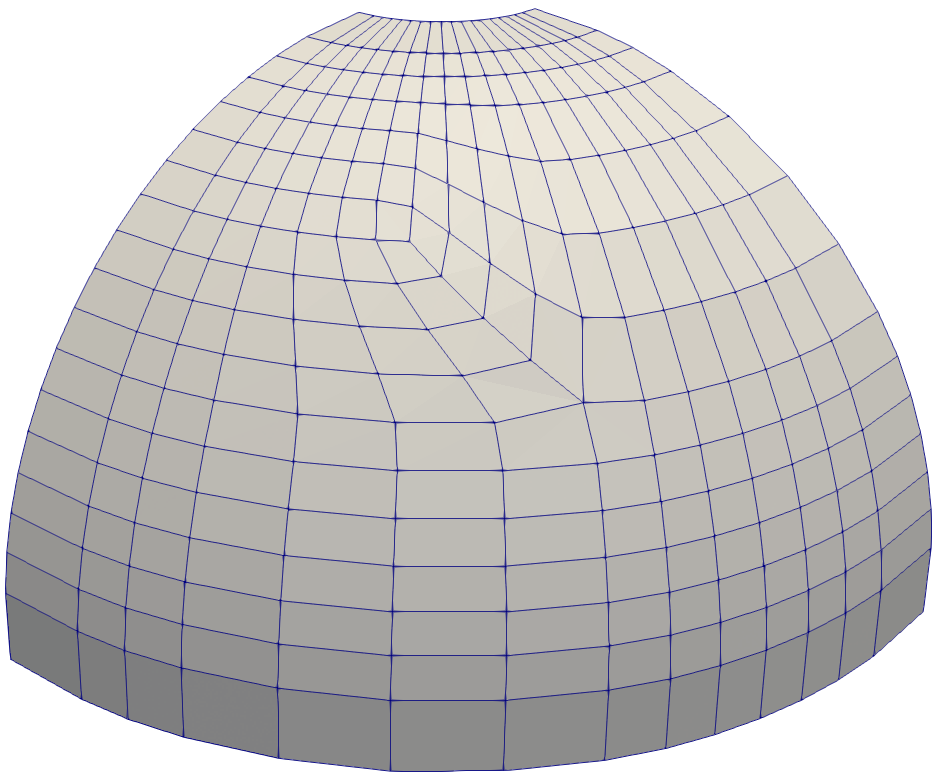}
    \caption{The initial control mesh of a hemisphere with two EPs }
    \label{fig:hemisphere_geo}
\end{figure}

We first study a 3D benchmark test to discuss the effect of nonlinearity on the proposed method. The geometry of interest is a hemisphere represented
 by ASUTS, where two EPs are involved; see Fig. \ref{fig:hemisphere_geo}. 
The time step is chosen as $\Delta t = 0.1$. The Dirichlet boundary \(\partial\Omega_D\) consists of three edges,  \(x=0\), \(y=0\), and \(z=0\), where we impose \(T=0\). We adopt the manufactured solution \(T(\mathbf{x}, t)=\left(1-e^{-5 t}\right) x y z\).  

We test three typical forms of thermal conductivity to cover linear and strongly nonlinear scenarios, where the constant conductivity \(k(T)=1\) serves as the linear baseline for comparison, the exponential nonlinear conductivity \(k(T)=e^{-T}\) represents smooth nonlinearity, and the periodic nonlinear conductivity \(k(T)=\sin(\pi T)\) corresponds to the strongly nonlinear scenario. The \(L^2\)-norm error of the steady-state temperature field and the total number of Newton iterations throughout the simulation are summarized in Table~\ref{tab:nonlinear_heat}, where the results of the proposed WQ are compared with those of the standard GQ. Note
 that each iteration needs one complete formation and assembly of the tangent stiffness matrix.

\begin{table}[htbp]
\centering
\caption{\(L^2\)-norm error and the total number of Newton iterations throughout the simulation using WQ and GQ  in the nonlinear heat transfer problem with different thermal conductivities \(k(T)\).}
\label{tab:nonlinear_heat}
\begin{tabular}{lcc}
\toprule
\(k(T)\) & $L^2$ error of WQ (\#iterations)& $L^2$ error  of GQ (\#iterations)\\
\midrule
\(1\) & 0.0185539 (50) & 0.0185539 (50) \\
\(e^{-T}\) & 0.0196745 (158) & 0.0196188 (158) \\
\(\sin(\pi T)\) & 0.0104808 (1600) & 0.0111105 (1532) \\
\bottomrule
\end{tabular}
\end{table}

We observe that the proposed WQ achieves nearly identical performance in terms of accuracy and convergence to the standard GQ in all cases of $k(T)$. This shows that the proposed method performs consistently well in practice. Moreover, the $9$ times speedup of matrix assembly observed in linear  problems is fully retained in nonlinear simulations. 

\begin{figure}[htbp]
\centering
\includegraphics[width=0.8\textwidth]{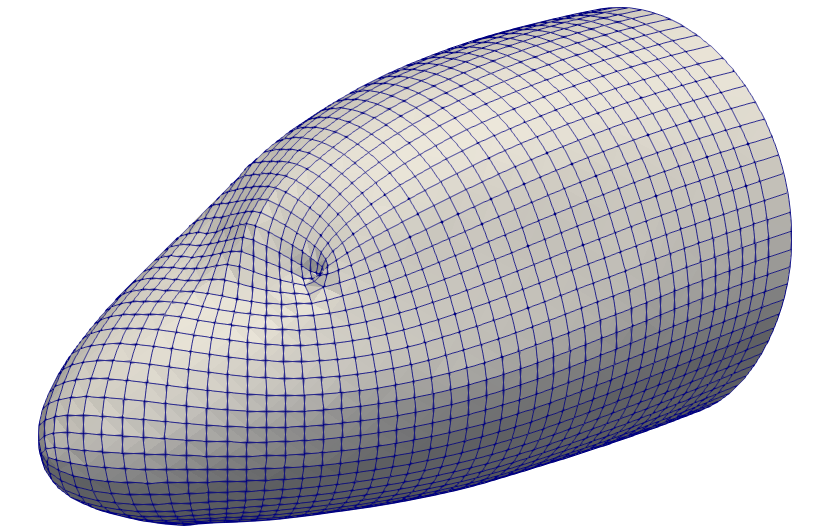}
\caption{The initial control mesh of an aircraft nose model, which involves clustered EPs in a single element.
}
\label{fig:aircraft_geo}
\end{figure}

As the last example, we apply the proposed  method to a  complex 3D surface to demonstrate its applicability for complex geometric models. The model is a  aircraft nose represented as a smooth ASUTS surface, where clustered EPs appear in a single element; see Fig. \ref{fig:aircraft_geo}.

We solve for the steady-state solution to the nonlinear heat transfer with a temperature-dependent thermal conductivity of aluminum alloy,  \(k(T)=200-0.05 T\).  A fixed temperature \(T=0\) is prescribed at the root of the nose surface. The source term has two parts. The convective part takes a simplified form,  \(f_{\text{conv}}=-3375-105T\)  \cite{Incropera2011}. The intensive aerodynamic heating flux \(f_{\text{aero}}=10^6\) is applied at the nose tip.  The time step is set to $\Delta t = 0.002$.

\begin{figure}[htbp]
\centering

  \begin{subfigure}{0.5\textwidth}
    \centering
    \includegraphics[width=\linewidth]{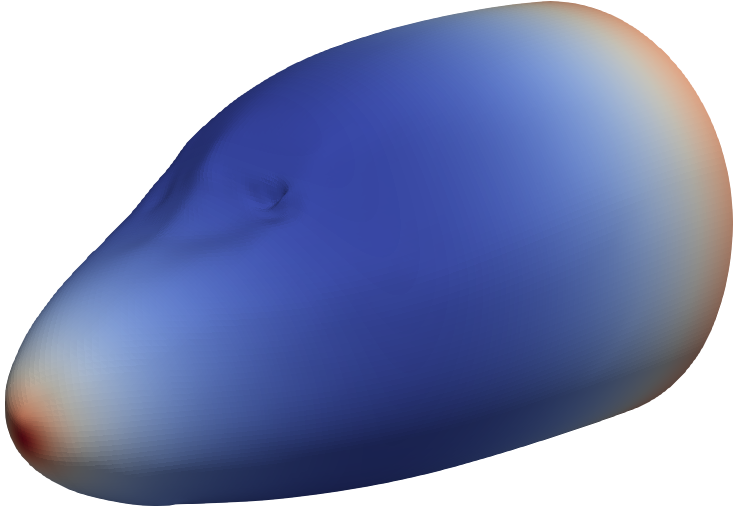}
    \caption{Proposed WQ}
  \end{subfigure}%
  \hfill
  \begin{subfigure}{0.5\textwidth}
    \centering
    \includegraphics[width=\linewidth]{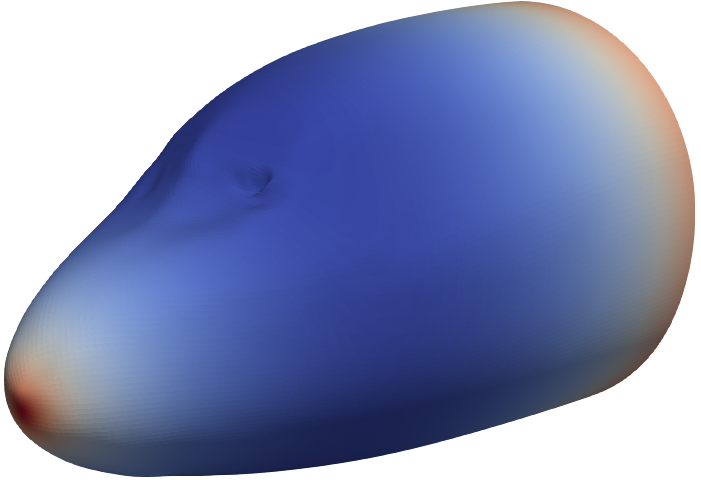}
    \caption{Standard GQ}
  \end{subfigure}

\hfill\\
\begin{minipage}[c]{0.6\textwidth}
  \centering
  \begin{subfigure}{\textwidth}
    \centering
    \includegraphics[width=\linewidth]{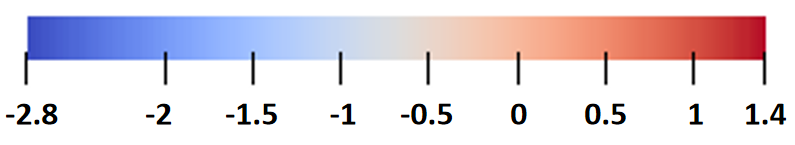}
    \caption{Color bar}
  \end{subfigure}
\end{minipage}

\caption{Steady-state temperature distributions using the proposed WQ (a) and the standard GQ (b).}
\label{fig:aircraft_temp}
\end{figure}

\begin{figure}[htbp]
\centering
\includegraphics[width=0.7\textwidth]{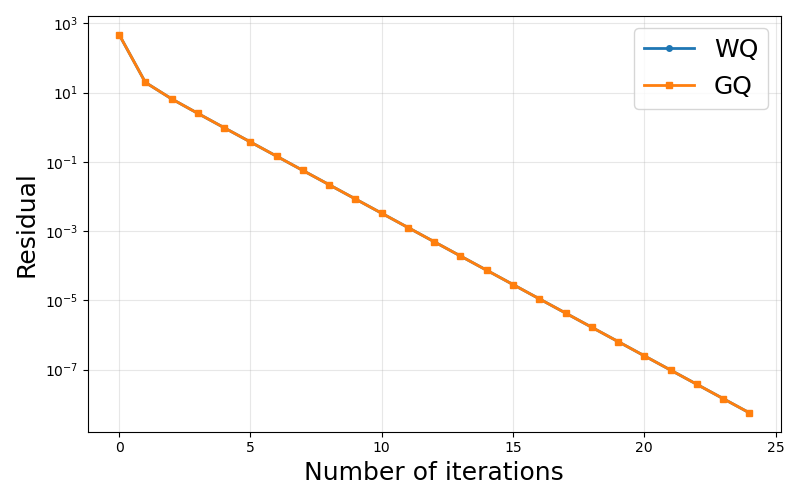}
\caption{Residual evolution with respect to the number of iterations using WQ and GQ.
}
\label{fig:aircraft_residual}
\end{figure}

The steady-state temperatures are computed by the proposed WQ  and the standard GQ. The results are shown in Fig. \ref{fig:aircraft_temp}. We observe that the two solutions are visually identical. The evolution of the  residual with respect to the number of iterations is plotted in Fig. \ref{fig:aircraft_residual}, where the results of  WQ and GQ  completely overlaps. It verifies that the proposed method can also deal with nonlinear problems defined on 3D complex surfaces without loss of accuracy.

\clearpage

\section{Conclusion and future work}
This work generalizes weighted quadrature to unstructured splines, facilitating fast formation and assembly of Galerkin matrices for complex geometries. To address the lack of a global parametric domain, we construct quadrature rules directly in the physical domain. We then introduce a general ``min-max" scheme to place quadrature points, ensuring the well-posedness condition with a minimal number of quadrature points. Given the underdetermined nature and the properties of unstructured splines, we further handle the ill-conditioning issue of the coefficient matrix through the truncated Singular Value Decomposition. In the end, the accuracy and efficiency of the proposed method are demonstrated in a variety of problems, including linear/nonlinear problems, high-order problems, and nonlinear problems on complex 3D surfaces. 

In the future, we plan to solve high-order nonlinear problems, such as the
 Cahn-Hilliard equation, based on the proposed method, where treatment of T-junctions will be investigated as well to support adaptivity. Moreover, applying it to topology optimization also shows great promise in gaining significant speedup.

\section*{Acknowledgements}
 Ji Sheng and Xiaodong Wei are partially supported by National Natural Science Foundation of China (No. 12494550/12494555 and No. 12571408).

\bibliography{cas-refs}
\end{document}